\documentclass[12pt]{amsart}

\usepackage{latexsym}
\usepackage{amssymb}
\usepackage[cp850]{inputenc}
\usepackage{epsfig}
\usepackage{psfrag}
\usepackage{amsthm}
\usepackage{amscd}
\usepackage{amsmath}
\usepackage{amsfonts}
\usepackage{graphics}
\usepackage[all]{xy}

\newtheorem{thm}{Theorem}[section]
\newtheorem{lem}{Lemma}[section]
\newtheorem{prop}{Proposition}[section]
\newtheorem{cor}{Corollary}[section]

\let\ssection=\section
\renewcommand{\section}{\setcounter{equation}{0}\ssection}

\newtheorem*{namedtheorem}{\theoremname}
\newcommand{\theoremname}{testing}

\theoremstyle{remark}

\newcommand{\BR}{\mathbb R}			
			
			\newcommand{\BZ}{\mathbb Z}

\DeclareMathOperator{\al}{\alpha}

\begin{document}

\title[Finiteness properties for a subgroup]{Finiteness properties for a subgroup of the pure symmetric automorphism group}
\author[A.~Pettet]{Alexandra Pettet}
\address{Dept.\ of Mathematics \\
  University of Michigan\\
  Ann Arbor, MI 48019}
\email{apettet@umich.edu}\begin{abstract}
Let $F_n$ be the free group on $n$ generators, and let ${\rm P}\Sigma_n$ be the group of automorphisms of $F_n$ that send each generator to a conjugate of itself. The kernel $K_n$ of the homomorphism ${\rm P}\Sigma_n \to {\rm P}\Sigma_{n-1}$, induced by mapping one of the free group generators to the identity, is finitely generated. We show that $K_n$ has cohomological dimension $n-1$, and that $H_i(K_n;\BZ)$ is not finitely generated for $2 \leq i \leq n-1$. It follows that $K_n$ is not finitely presentable for $n \geq 3$.
\end{abstract}
\thanks{\tiny The author is partially supported by NSF grant
  DMS-0856143 and NSF RTG DMS-0602191}

\date{\today}

\maketitle

\section{Introduction}
\label{}

In recent work, Brendle-Hatcher \cite{Brendle-Hatcher} proved that the space of all smooth links in $\BR^3$ isotopic to the trivial link of $n$ components has the homotopy type of the finite dimensional subspace of configurations of $n$ unlinked circles, and thus their fundamental groups are isomorphic. The fundamental group of the latter space is a 3-dimensional analogue of the classical braid group (the space of configurations of $n$ points in $\BR^2$), and Goldsmith \cite{Goldsmith} showed that it is isomorphic to the {\it symmetric automorphism group}, the group of automorphisms of $F_n$ which send every generator to a conjugate of another generator or its inverse.

The subgroup consisting of those automorphisms which send every generator to a conjugate of itself (or, in mapping class group terms, those classes which send every oriented circle in $\BR^3$ back to itself) is known as the {\it pure symmetric automorphism group}, denoted by ${\rm P}\Sigma_n$. McCool \cite{McCool} gave a finite presentation for ${\rm P}\Sigma_n$, and Brownstein-Lee \cite{Brownstein-Lee} computed its cohomology when $n=3$. Collins \cite{Collins} proved that ${\rm P}\Sigma_n$ has cohomological dimension $n-1$; it also follows from his work that ${\rm P}\Sigma_n$ is $FP_\infty$.
Later, Brady-McCammond-Meier-Miller \cite{Brady} showed that ${\rm P}\Sigma_n$ is a duality group, and Jensen-McCammond-Meier \cite{Jensen-McCammond-Meier} determined completely the structure of the cohomology ring of ${\rm P}\Sigma_n$ for $n \geq 3$. 

Let $PB_n$ denote the pure braid group, the elements of the braid group that send each puncture back to itself. It is well-known that for all $n$ there is a homomorphism $\pi:PB_n\to PB_{n-1}$ induced by ``filling in'' a puncture. In fact, there is the following split exact sequence: 
\begin{equation}\label{eq1}
\xymatrix{ 1 \ar[r] & F_{n-1} \ar[r] & PB_n \ar[r]^{\pi} & PB_{n-1} \ar[r] & 1} 
\end{equation}
In particular, the pure braid group may be regarded as an iteration of semi-direct products of free groups. The pure braid group $PB_n$ is isomorphic to a subgroup of ${\rm P}\Sigma_n$, and by ``filling in'' the $n$th circle we obtain a split exact sequence compatible with (\ref{eq1}): 
\begin{equation}\label{eq2}
\xymatrix{1 \ar[r] & K_n \ar[r] & {\rm P}\Sigma_n \ar[r]^{\pi} & {\rm P}\Sigma_{n-1} \ar[r] & 1}
\end{equation}
For $n=2$ the kernel $K_2$ is equal to ${\rm P}\Sigma_2$. For $n \geq 2$, the group $K_n$ is finitely generated (compare with Lemma \ref{H1} below), and hence $H_1(K_n;\BZ)$ is finitely generated. The main purpose of this note is to study the higher homology groups of $K_n$ for $n \geq 3$: 

\begin{thm}\label{main}
The group $K_n$ has cohomological dimension $n-1$. For $n \geq 3$ its $i$th homology group $H_i(K_n;\BZ)$ is not finitely generated for $2 \leq i \leq n-1$.
\end{thm}
\noindent Collins-Gilbert proved that $K_3$ is not finitely presentable in \cite{Collins-Gilbert}. Theorem \ref{main} yields an independent proof of this fact, generalizing to all $n \geq 3$: 
\begin{cor} 
$K_n$ is not finitely presentable for $n \geq 3$. 
\end{cor}
\noindent As pointed out by Brendle-Hatcher \cite{Brendle-Hatcher}, the corollary suggests that these kernels $K_n$ are unlikely to have nice interpretations in terms of configuration spaces of circles.

\medskip
The author is grateful to her advisor Benson Farb for his unfailing enthusiasm and guidance, and to Allen Hatcher, John Meier, and Juan Souto for helpful suggestions. The author was partially supported by the Natural Science and Engineering Council of Canada (NSERC) Post Graduate Scholarship during the conception of this note. 

\section{Finitely generated homology groups}\label{fg-homology}%
In this section we verify the finite generation of $K_n$ and compute its first homology group, $H_1(K_n;\BZ)$, and its cohomological dimension.
\begin{lem}\label{H1}
The group $K_n$ is finitely generated, and its first homology group is 
\[ H_1(K_n;\BZ) \simeq \BZ^{2n-2} \]
\end{lem}
\begin{proof} McCool \cite{McCool} proved that the group ${\rm P}\Sigma_n$ is generated by 
\[ 
\al_{ij}(x_r) = 
\left \{
\begin{array}{ll}
x_r & r \neq i \\*
x_j x_i x_j^{-1} & r = i
\end{array}
\right.
\]
with relators
\[ [\al_{ij},\al_{kl}], \quad [\al_{ik},\al_{jk}], \quad [\al_{ij},\al_{ik}\al_{jk}] \]
for distinct $i,j,k$, and $l$. 
\noindent It is clear that $K_n$ is normally generated by 
\[ \{\al_{in}, \al_{ni} \ | \ 1 \leq i \leq n-1\} \] 
In fact by examining the McCool generators, these elements are seen to generate $K_n$: 
\begin{eqnarray*} & \al_{ij}^{\pm} \al_{ni}\al_{ij}^{\mp} = \al_{nj}^{\mp}\al_{ni}\al_{nj}^{\pm}  \quad \quad \al_{jk}^{\pm}\al_{ni}\al_{jk}^{\mp} = \al_{ni} \quad \quad \al_{ji}^{\pm} 			\al_{ni} \al_{ji}^{\mp} = \al_{ni} & \\
& \al_{ij}^{\pm} \al_{in}\al_{ij}^{\mp} = \al_{nj}^{\mp}\al_{in}\al_{nj}^{\pm} \quad \quad \al_{jk}^{\pm}\al_{in}\al_{jk}^{\mp} = \al_{in} & \\
& \al_{ji}^{-1} \al_{in} \al_{ji} = \al_{ji}^{-1} \al_{jn}^{-1} \al_{ji} \al_{in} \al_{jn} 
						=  \al_{ni} \al_{jn}^{-1} \al_{ni}^{-1} \al_{in} \al_{jn}  &\\						
& \al_{ji} \al_{in} \al_{ji}^{-1} = \al_{jn} \al_{in} \al_{ji} \al_{jn}^{-1} \al_{ji}^{-1} 
						= \al_{jn} \al_{in} \al_{ni}^{-1} \al_{jn}^{-1} \al_{ni} & \end{eqnarray*}
The last two expressions are each derived from a conjugate of a McCool relator, followed by a substitution using a second relator. 

Consider the free group $F(\{ \al_{in}, \al_{ni} \})$ of rank $2n-2$ on the generators of $K_n$, a subgroup of the free group $F(\{\al_{ij}\})$ of rank $n^2-n$ on the generators of ${\rm P}\Sigma_n$. It is clear from McCool's presentation that the kernel of the map $F(\{\al_{ij}\}) \to {\rm P}\Sigma_n$ is contained in the commutator subgroup. An element in the kernel of $F(\{\al_{in}, \al_{ni} \}) \to K_n$ lies in 
\[ [F(\{\al_{ij}\}), F(\{\al_{ij}\})] \cap F(\{ \al_{in}, \al_{ni} \}) \]
and such an element must also lie in the commutator subgroup of $F(\{ \al_{in}, \al_{ni} \})$. This shows that $H_1(K_n;\BZ) \simeq H_1(F(\{ \al_{in}, \al_{ni} \});\BZ)$, completing the proof of Lemma \ref{H1}. 
\end{proof}

Jensen-Wahl \cite{Jensen-Wahl} describe an $(n-1)$-dimensional contractible simplicial complex $X_n$ on which ${\rm P}\Sigma_n$ acts freely with compact quotient. Briefly, this complex $X_n$ is the geometric realization of the poset of symmetric based graphs with fundamental group $F_n$, and a marking from a basis $\{ x_1, \ldots, x_n \}$ to each graph $\Gamma$ which induces an isomorphism $F_n \to \pi_1(\Gamma)$. A symmetric graph is one in which every edge belongs to a unique cycle, and the partial ordering is given by the collapsing of edges. The complex $X_n$ embeds into the spine of Autre Space, the based-graph version of Culler-Vogtmann's Outer space. We refer the reader to \cite{Jensen-Wahl} for details. 
\begin{lem}\label{highercohom}
 $K_n$ has cohomological dimension $n-1$.
\end{lem}
\begin{proof}
The existence of $X_n$ gives an upper bound of $n-1$ for the cohomological dimension of $K_n$. The elements  $\{\al_{1n}, \ldots \al_{n-1,n}\}$ generate a free abelian subgroup of rank $n-1$, so that $n-1$ is also a lower bound. This completes the proof of the lemma, and thereby the first part of Theorem \ref{main}.
\end{proof}

\section{Non-finitely generated homology groups}
We begin with a short lemma about $H_{n-1}(K_n;\BZ)$:

\begin{lem}\label{nontrivial}
The group $H_{n-1}(K_n;\BZ)$ has a nontrivial element.
\end{lem}
\begin{proof}
The subgroup $K_n$ containes the $n-1$ commuting elements $\al_{1n}, \ldots, \al_{n-1,n}$. From the McCool relations, we can verify that we have homomorphisms
\begin{eqnarray*}
\xymatrix{\BZ^{n-1} \ar[r] & K_n \ar[r] & \BZ^{n-1}}
\end{eqnarray*}
whose composition is the identity. Therefore the induced map $H_{n-1}(\BZ^{n-1};\BZ) \to H_{n-1}(K_n;\BZ)$ is injective.
\end{proof}

We next prove a proposition which, together with Lemma \ref{nontrivial}, proves the theorem. The author is thankful to A.~Hatcher for suggesting this proposition as a very nice simplification of arguments in an earlier version of this note: 

\begin{prop}\label{keyprop} Let $\Gamma$ be a group acting freely and simplicially on a contractible $(n-1)$-dimensional complex $X$, and let $K$ be normal subgroup of $\Gamma$ of infinite index. Then if $H_{n-1}(K;\BZ)$ is nonzero, it is not finitely generated.
\end{prop}
\begin{proof}
By assumption, $K$ acts freely on the contractible complex $X$, so $Y=X/K$ is an Eilenberg-MacLane space of type $K(K,1)$. Thus by the assumption that $H_{n-1}(K;\BZ) \neq 0$, we have $H_{n-1}(Y;\BZ) \neq 0$. A nontrivial $(n-1)$-cycle of $Y$ is represented by a finite sum of $(n-1)$-simplices, so there exists a nontrivial finite subcomplex $A$ of $Y$ such that $H_{n-1}(A;\BZ) \neq 0$. As $\Gamma/K$ acts freely on $Y$, and as $K$ has infinite index in $\Gamma$, there is an infinite set of pairwise disjoint translates of $A$ by $\Gamma/K$; denote the union of such a set of translates by $U$. Clearly $H_{n-1}(U;\BZ)$ is not finitely generated. The proof is complete by the following exact sequence on the relative pair $(Y,U)$:
\begin{eqnarray*}
\xymatrix{ \cdots \ar[r] & H_n(Y,U;\BZ) \ar[r] & H_{n-1}(U;\BZ) \ar[r] & H_{n-1}(Y;\BZ) \ar[r] & H_{n-1}(Y,U;\BZ) \ar[r] & \cdots} 
\end{eqnarray*}
The first term $H_n(Y,U;\BZ)=0$ as $Y$ has dimension $n-1$.
\end{proof}
For $n \geq 3$, the subgroup $K_n$ has infinite index in ${\rm P}\Sigma_n$, and so Lemma \ref{nontrivial} and Proposition \ref{keyprop} applied to the Jensen-Wahl complex $X_n$ imply that $H_{n-1}(K_n;\BZ)$ is not finitely generated. Now observe that there exists a split homomorphism $K_n \to K_{n-1}$. For all $i$, this induces maps $H_i(K_{n-1}; \BZ) \to H_i(K_n; \BZ)$ and $H_i(K_n; \BZ) \to H_i(K_{n-1}; \BZ)$ whose composition $H_i(K_{n-1};\BZ) \to H_i(K_{n-1};\BZ)$ is the identity; the first map must be injective. Theorem \ref{main} then holds for $n\geq 3$ by induction on $n$.


\end{document}